\documentclass[11pt]{article}
\usepackage{mathrsfs,amsmath,amssymb,epsfig,subfigure,color}

\renewcommand{\arraystretch}{1.2}

\newcommand{\BE}{\begin{equation}}
\newcommand{\EE}{\end{equation}}

\newtheorem{theorem}{Theorem}[section]
\newtheorem{lemma}{Lemma}[section]
\newtheorem{example}{Example}[section]

\newtheorem{remark}{Remark} [section]

\def\dfrac{\displaystyle\frac}
\def\dsum{\displaystyle\sum}
\def\dint{\displaystyle\int}

\begin{document}
\title{\textbf{Compact difference schemes for the modified anomalous fractional sub-diffusion equation
and the fractional diffusion-wave equation}}
\author{Zhibo Wang\thanks{Corresponding author. Email: zhibowangok@gmail.com. Department of Mathematics, University of Macau, Av. Padre Tom\'{a}s Pereira Taipa, Macau, China}
\and Seakweng Vong\thanks{ Email: swvong@umac.mo. Department of Mathematics, University of Macau, Av. Padre Tom\'{a}s Pereira Taipa, Macau, China}}

 \date{}
 \maketitle
 \begin{abstract}
 In this paper, compact finite difference schemes for the
modified anomalous fractional sub-diffusion equation and fractional
diffusion-wave equation are studied. Schemes proposed previously can
at most achieve temporal accuracy of order which depends on the
order of  fractional derivatives in the equations and is usually
less than two. Based on the idea of weighted and shifted
Gr\"{u}nwald difference operator, we establish schemes with
temporal and spatial accuracy order equal to two and four
respectively.
\end{abstract}
 {\bf Keywords:} Modified anomalous fractional sub-diffusion equation, Fractional diffusion-wave equation, Compact difference scheme, Weighted and shifted Gr\"{u}nwald difference operator
\section{Introduction}
 Since fractional differential equations turn out to model many
 physical processes more accurately than the classical ones, in the past decades,
 increasing attentions  have  been made on these equations.
 Readers can refer to the books \cite{Podlubny,Kilbas} for
 theoretical results on fractional differential equations.
 This paper concerns with methods for obtaining accurate numerical
 approximations to the solutions of fractional sub-diffusion equations and fractional diffusion-wave  equations.
 A fractional sub-diffusion equation is an  integro-partial differential equation obtained from the classical diffusion
 equation by replacing the first-order time derivative by a fractional
 derivative of order between zero and one. When the time derivative
 is of order between one and two, we get a fractional diffusion-wave equation. Fractional
 derivatives of order between zero and one are widely used in describing anomalous diffusion
 processes \cite{Klages}, while fractional diffusion-wave equations have applications in modeling  universal electromagnetic,
 acoustic, and mechanical responses \cite{Nigmatullin1, Nigmatullin2}.

 Numerical methods for the modified anomalous fractional sub-diffusion equation and the diffusion-wave equation have been considered by many authors, one
 may refer to \cite{ChenC}--\cite{r6} and the references therein.
 We point out here that one of the main tasks for developing accurate finite difference scheme of fractional
 differential equations is to discretize the fractional derivatives.
 Noticing that fractional derivatives are defined through integrals, one can
 approximate the derivatives by interpolating polynomials
 \cite{Langlands}. Using this idea, Sousa and
 Li developed a second order discretization for the Riemann-Liouville
 fractional derivative \cite{Sousaa}.

 We note that there are alternative ways to tackle the problem. In \cite{Meerschaert},
  a shifted Gr\"{u}nwald formula was
 proposed by Meerschaert and Tadjeran to approximate fractional
 derivatives of order $\alpha\in(1,2)$ for fractional
 advection-dispersion flow equations. It is also worth to mention
 that stability of forward-Euler scheme and weighted averaged difference scheme based on
 Gr\"{u}nwald-Letnikov approximation were analyzed in \cite{Yuste1, Yuste2} for time fractional diffusion equations.
 Very recently, accurate finite difference schemes based on  weighted and shifted Gr\"{u}nwald
 difference operator were developed for solving space fractional diffusion
 equations in \cite{Tian, Zhou}.

 Inspired by their work on the weighted and shifted Gr\"{u}nwald
 difference operator, in this paper, we establish high order schemes for
 the fractional diffusion-wave equation and the modified anomalous fractional
sub-diffusion equation, which were proposed and studied in
\cite{SunWu, Du} and \cite{F-Liu}--\cite{Mohebbi}, respectively.
 We remark that the order of temporal accuracy of schemes
 proposed previously can at most be a fraction depending on the fractional derivatives in the
 equations and is usually less than two.
 We show that the  schemes proposed in this paper are of order $\tau^2+h^4$,
 where $\tau$ and $h$ are the temporal and  spatial step sizes respectively.

 This paper is organized as follows. Some preliminaries will be
 given in the next section. Compact schemes are proposed and studied in
 Section 3 and 4 for the modified anomalous fractional
 sub-diffusion equation and the fractional diffusion-wave equation
 respectively. In the last section, numerical tests are carried out
 to justify the theoretical results.

 \section{Preliminaries}
 We first recall that the Caputo fractional derivative of order $\gamma\in(1,2)$ for a function $f(t)$  is defined as
 $$_{a}^CD_t^{\gamma}f(t)=
 \frac1{\Gamma(2-\gamma)}\int_a^t
 \frac{f''(s)}{(t-s)^{\gamma-1}}ds,$$
 with $\Gamma(\cdot)$ being the gamma function and,
 for $\alpha\in(0,1)$, the Riemann-Liouville fractional
 derivative of order $\alpha$ for $f$ is defined as
 $$_{a}D_t^{\alpha} f(t)=
 \frac1{\Gamma(1-\alpha)}\dfrac{d}{dt}\int_a^t
 \frac{f(s)}{(t-s)^{\alpha}}ds.$$
 Closely related to the fractional derivatives of a function
 is the Riemann-Liouville fractional integral which is given by
 $$_{a}I_t^\alpha f(t)=\dfrac1{\Gamma(\alpha)}\dint_a^t\frac{f(s)}{(t-s)^{1-\alpha}}ds.$$
 In order to develop second order approximation  of the Riemann-Liouville fractional derivative, we consider
 the shifted Gr\"{u}nwald approximation \cite{Meerschaert} to the Riemann-Liouville fractional derivative given by
 $$\mathcal{A}_{\tau,r}^{\alpha} f(t)=\tau^{-\alpha}\dsum_{k=0}^{\infty}g_k^{(\alpha)}f(t-(k-r)\tau),$$
 where $g_k^{(\alpha)}=(-1)^k\binom{\alpha}{k}$ for $k\geq0$.
 Inspired by \cite{Lubich}, we similarly introduce the shifted operator to the Riemann-Liouville fractional integral defined by
 $$\mathcal{B}_{\tau,r}^\alpha f(t)=\tau^\alpha\dsum_{k=0}^{\infty}\omega_k^{(\alpha)}f(t-(k-r)\tau),$$
 where $\omega_k^{(\alpha)}=(-1)^k\binom{-\alpha}{k}$
 for $k\geq0$.

 The following lemma was given in {\rm\cite{Ervin}}.
 \begin{lemma}\label{fourior}
 Suppose $\alpha>0,~f(t)\in L^p(\mathbb{R}),~p\geq1$. The Fourier transform of the Riemann-Liouville fractional integral satisfy the following:
 $$\mathscr{F}[{_{-\infty}I_t^\alpha f(t)}]=(i\omega)^{-\alpha}\hat{f}(\omega),$$
 where $\hat{f}(\omega)=\int_\mathbb{R}e^{-i\omega t}f(t)dt$ denotes the Fourier transform of $f(t)$.
 \end{lemma}
 Lemma \ref{fourior} can be used to obtain

 \begin{lemma}\label{main-operator}
 {\rm(i)} Let $f(t)\in L^1(\mathbb{R}),~_{-\infty}D_t^{\alpha+2}f$
 and its Fourier transform belong to $L^1(\mathbb{R})$, and define the weighted and shifted Gr\"{u}nwald difference operator by
 \begin{equation}\label{Dengoperator}
 \mathcal{D}_{\tau,p,q}^{\alpha} f(t)
 =\dfrac{\alpha-2q}{2(p-q)}\mathcal{A}_{\tau,p}^{\alpha} f(t)
 +\dfrac{2p-\alpha}{2(p-q)}\mathcal{A}_{\tau,q}^{\alpha} f(t),
 \end{equation}
 then we have
 $$\mathcal{D}_{\tau,p,q}^{\alpha} f(t)={_{-\infty}D_t^\alpha f(t)}+O(\tau^2)$$
 for $t\in \mathbb{R}$, where $p$ and $q$ are integers and $p\neq q$.

 {\rm(ii)} Let $f(t),~{_{-\infty}I_t^\alpha f(t)}$ and
 $(i\omega)^{2-\alpha}\mathscr{F}[f](\omega)$
 belong to $L^1(\mathbb{R})$. Define the weighted and shifted difference operator by
 $$\mathcal{I}_{\tau,p,q}^\alpha f(t)
 =\dfrac{2q+\alpha}{2(q-p)}\mathcal{B}_{\tau,p}^\alpha f(t)
 +\dfrac{2p+\alpha}{2(p-q)}\mathcal{B}_{\tau,q}^\alpha f(t),$$
 then we have
 $$\mathcal{I}_{\tau,p,q}^\alpha f(t)={_{-\infty}I_t^\alpha f(t)}+O(\tau^2)$$
 for $t\in \mathbb{R}$, where $p$ and $q$ are integers and $p\neq q$.
 \end{lemma}

 {\bf Proof.} The proof of (i) can be found in \cite{Tian}. The
 proof of (ii) is similar to that of (i) but we give it here for the
 completeness of our presentation.

 Referring to the definition of $\mathcal{B}_{\tau,r}^\alpha$, we let
 \begin{equation}\label{weight}
 \mathcal{I}_{\tau,p,q}^\alpha f(t)=\tau^\alpha\bigg[\mu_1\dsum_{k=0}^{\infty}\omega_k^{(\alpha)}f(t-(k-p)\tau)
 +\mu_2\dsum_{k=0}^{\infty}\omega_k^{(\alpha)}f(t-(k-q)\tau)\bigg].
 \end{equation}
 Taking Fourior transform on (\ref{weight}), we get
 \begin{equation}\label{fourior-weight}
 \begin{array}{rl}
 \mathscr{F}[\mathcal{I}_{\tau,p,q}^\alpha f](\omega)
 &=\tau^\alpha\dsum_{k=0}^{\infty}\omega_k^{(\alpha)}\big[\mu_1e^{-i\omega(k-p)\tau} +\mu_2e^{-i\omega(k-q)\tau}\big]\mathscr{F}[f](\omega)\\
 &=\tau^\alpha\big[\mu_1(1-e^{-i\omega\tau})^{-\alpha}e^{i\omega\tau p}+\mu_2(1-e^{-i\omega\tau})^{-\alpha}e^{i\omega\tau q}\big]\mathscr{F}[f](\omega)\\
 &=(i\omega)^{-\alpha}[\mu_1W_p(i\omega\tau)+\mu_2W_q(i\omega\tau)]\mathscr{F}[f](\omega),
 \end{array}
 \end{equation}
 where
 \begin{equation}\label{taylor}
 W_r(z)=z^\alpha(1-e^{-z})^{-\alpha}e^{rz}=1+\big(r+\frac\alpha2\big)z+O(z^2),
 \quad r=p,q.
 \end{equation}
 In order to achieve second order accuracy,
 we let the coefficients $\mu_1$ and $\mu_2$ satisfy the following system:
 $$\left\{
 \begin{array}{l}
 \mu_1+\mu_2=1, \\
 (p+\frac\alpha2)\mu_1+(q+\frac\alpha2)\mu_2=0,
 \end{array}\right.$$
 which implies that $\mu_1=\frac{2q+\alpha}{2(q-p)}$ and
 $\mu_2=\frac{2p+\alpha}{2(p-q)}$.

 Denote
 $\hat{g}(\omega,\tau)=\mathscr{F}[\mathcal{I}_{\tau,p,q}^\alpha f](\omega)-\mathscr{F}[_{-\infty}I_t^\alpha f](\omega)$,
 then by (\ref{fourior-weight}), (\ref{taylor}) and Lemma \ref{fourior}, we have
 $$|\hat{g}(\omega,\tau)|\leq
 C\tau^2|i\omega|^{2-\alpha}|\mathscr{F}[f](\omega)|.$$
 Thus
 $$|\mathcal{I}_{\tau,p,q}^\alpha f-{_{-\infty}I_t^\alpha f}|
 =|g|\leq\dfrac1{2\pi}\dint_\mathbb{R}|\hat{g}(\omega,\tau)|d\omega
 \leq C\|(i\omega)^{2-\alpha}\mathscr{F}[f](\omega)\|_{L^1}
 \tau^2=O(\tau^2).\qquad\Box$$

 We are now ready to establish our high order compact schemes in the
 next two sections.

 \section{A compact scheme for the fractional sub-diffusion
 equation}
 Consider the following modified anomalous fractional sub-diffusion equation
 \begin{equation}\label{sub-main1}
 \dfrac{\partial u(x,t)}{\partial t}
 =\bigg(\kappa_1\dfrac{\partial^{\alpha}}{\partial t^{\alpha}}+\kappa_2\dfrac{\partial^{\beta}}{\partial t^{\beta}}\bigg)\bigg[\dfrac{\partial^2u(x,t)}{\partial x^2}\bigg]+f(x,t),
  \quad 0\leq x\leq L, \quad 0<t\leq T,
 \end{equation}
 subject to
 \begin{equation}\label{sub-main2}
 \begin{array}{rl}
 &u(x,0)=0, \quad 0\leq x\leq L,\\
 &u(0,t)=\varphi_1(t), \quad u(L,t)=\varphi_2(t), \quad 0<t\leq T,
 \end{array}
 \end{equation}
 where $0<\alpha, \beta<1$, $\kappa_1, \kappa_2\geq0$.
 We have used $\frac{\partial^{\alpha}}{\partial t^{\alpha}}$ and $\frac{\partial^{\beta}}{\partial t^{\beta}}$ to
 denote the Riemann-Liouville fractional operators $_{0}D_t^{\alpha}$ and $_{0}D_t^{\beta}$ with respect to the time variable $t$.

 \begin{remark}
 {\rm(i)} We note that, in {\rm\cite{F-Liu, Q-Liu, Mohebbi}}, the fractional derivatives in the equation are of order $1-\tilde\alpha$ and $1-\tilde\beta$ with $\tilde\alpha,\tilde\beta\in(0,1)$. We have changed the notations here in order to match the presentation for the two types of equations discussed in this paper.

 {\rm(ii)} Without loss of generality, we have assumed the initial condition  $u(x,0)=0$. If $u(x,0)=\psi(x)$, one may consider the equation for $v(x,t)\dot=u(x,t)-\psi(x)$  instead.
 \end{remark}

 To develop a finite difference scheme for the problem (\ref{sub-main1})--(\ref{sub-main2}), we let $h=\frac LM$ and $\tau=\frac TN$ be the spatial and temporal step sizes respectively, where $M$ and $N$ are some given positive integers.
 For $i=0,1,\ldots,M$ and $k=0,1,\ldots,N$, denote $x_i=ih,~t_k=k\tau$.
 For any grid function $u=\{u_i^k|0\leq i\leq M,~ 0\leq k\leq N\}$, we introduce the following
 notations:
 $$\delta_xu_{i-\frac12}^k=\dfrac1h(u_i^k-u_{i-1}^k),~~\delta_x^2u_{i}^k=\dfrac1h(\delta_xu_{i+\frac12}^k-\delta_xu_{i-\frac12}^k),$$
 $${\cal H}u_i=\left\{
 \begin{array}{ll}
 \Big(1+\dfrac{h^2}{12}\delta_x^2\Big)u_i=\dfrac1{12}(u_{i-1}+10u_i+u_{i+1}), &1\leq i\leq M-1,\\
 u_{i}, &i=0 \mbox{ or } M.
 \end{array}\right.$$
 With this notations, we study the problem under the following inner product and norms:
 $$\langle u,v\rangle=h\dsum_{i=1}^{M-1}u_iv_i,
 ~~\langle\delta_xu,\delta_xv\rangle=h\dsum_{i=0}^{M-1}\big(\delta_xu_{i+\frac12}\big)\big(\delta_xv_{i+\frac12}\big),$$
 $$\|u\|^2=\langle u,u\rangle,
 ~~\|u\|_\infty=\max_{0\leq i\leq M}|u_i|.$$
 It is easy to check that, if $v\in{\cal V}=\{w|w=(w_0,w_1,\ldots,w_M),w_0=w_M=0\}$, the following identity holds
 $$\langle\delta_x^2u,v\rangle=-\langle\delta_xu,\delta_xv\rangle,$$
 which plays important role in our analysis.

 \bigskip
 Note that one can  continuously extend the solution $u(x,t)$  to be zero for $t<0$.
 By choosing $(p,q)=(0,-1)$ in (i) of Lemma \ref{main-operator},
 we get $\frac{\alpha-2q}{2(p-q)}=\frac{2+\alpha}2,~
 \frac{2p-\alpha}{2(p-q)}=-\frac{\alpha}2$ in (\ref{Dengoperator}), which gives
 $$\begin{array}{rl}
 \frac{\partial^{\alpha}}{\partial t^{\alpha}}[u_{xx}(x_i,t_{n+1})]
 &=\tau^{-\alpha}\bigg(\dfrac{2+\alpha}2\dsum_{k=0}^{n+1}g_k^{(\alpha)}\delta^2_xu^{n+1-k}_i
 -\dfrac{\alpha}2\dsum_{k=0}^{n}g_k^{(\alpha)}\delta^2_xu^{n-k}_i\bigg)+O(\tau^2+h^2)\\
 &=\tau^{-\alpha}\dsum_{k=0}^{n+1}\lambda_k^{(\alpha)}\delta^2_xu^{n+1-k}_i+O(\tau^2+h^2),
 \end{array}$$
 where
 \begin{equation}\label{sequence-lambda}
 \lambda_0^{(\alpha)}=\dfrac{2+\alpha}2g_0^{(\alpha)},
 ~~\lambda_k^{(\alpha)}=\dfrac{2+\alpha}2g_{k}^{(\alpha)}
 -\dfrac{\alpha}2g_{k-1}^{(\alpha)},~~k\geq1.
 \end{equation}

 We can therefore consider an weighted Crank-Nicolson type discretization for equation (\ref{sub-main1}) given by
 $$\begin{array}{ll}
 \dfrac{u^{n+1}_i-u^n_i}{\tau}=&\dfrac{\kappa_1\tau^{-\alpha}}2
 \bigg(\dsum_{k=0}^{n+1}\lambda_k^{(\alpha)}\delta^2_xu^{n+1-k}_i
 +\dsum_{k=0}^{n}\lambda_k^{(\alpha)}\delta^2_xu^{n-k}_i\bigg)\\
 &+\dfrac{\kappa_2\tau^{-\beta}}2
 \bigg(\dsum_{k=0}^{n+1}\lambda_k^{(\beta)}\delta^2_xu^{n+1-k}_i
 +\dsum_{k=0}^{n}\lambda_k^{(\beta)}\delta^2_xu^{n-k}_i\bigg)
 +\dfrac12(f^n_i+f^{n+1}_i).
 \end{array}$$

 \smallskip
 In order to raise the accuracy in the spatial direction, we need the following lemma:
 \begin{lemma}\label{compact}{\rm(\cite{Gao})}
 Denote $\zeta(s)=(1-s)^3[5-3(1-s)^2]$. If $f(x)\in {\cal C}^6[x_{i-1},x_{i+1}],~1\leq i\leq M-1$, then it holds that
 $$\begin{array}{rl}
 \dfrac1{12}[f''(x_{i-1})+10f''(x_i)+f''(x_{i+1})]
 =&\dfrac1{h^2}[f(x_{i-1})-2f(x_i)+f(x_{i+1})]\\
 &+\dfrac{h^4}{360}\dint_0^1[f^{(6)}(x_i-sh)+f^{(6)}(x_i+sh)]\zeta(s)ds.
 \end{array}$$
 \end{lemma}

 Based on  Lemma \ref{compact}, we therefore propose the following compact scheme:
 \begin{align}\label{compact-scheme1}\nonumber
 &{\cal H}(u^{n+1}_i-u^n_i)=\dfrac{\kappa_1\tau^{1-\alpha}}2
 \bigg(\dsum_{k=0}^{n+1}\lambda_k^{(\alpha)}\delta^2_xu^{n+1-k}_i
 +\dsum_{k=0}^{n}\lambda_k^{(\alpha)}\delta^2_xu^{n-k}_i\bigg)\\\nonumber
 &\quad+\dfrac{\kappa_2\tau^{1-\beta}}2
 \bigg(\dsum_{k=0}^{n+1}\lambda_k^{(\beta)}\delta^2_xu^{n+1-k}_i
 +\dsum_{k=0}^{n}\lambda_k^{(\beta)}\delta^2_xu^{n-k}_i\bigg)
 +\dfrac\tau2{\cal H}(f^n_i+f^{n+1}_i),\\
 &0\leq n\leq N-1,~~1\leq i\leq M-1,\\\nonumber
 &u_0^n=\varphi_1^n,\quad u_M^n=\varphi_2^n,\quad 1\leq n\leq N,\\\label{compact-scheme2}
 &u_i^0=0,\quad 0\leq i\leq M.
 \end{align}

 It is easy to see that at each time level, the difference scheme (\ref{compact-scheme1})--(\ref{compact-scheme2}) is a linear tridiagonal system with strictly diagonal dominant coefficient matrix, thus the difference scheme has a unique solution.

 The following lemmas are critical for establishing the convergence of the proposed scheme.
 \begin{lemma}\label{sub-sequence}
 Let $\big\{\lambda_n^{(\alpha)}\big\}_{n=0}^\infty$ be defined as in {\rm(\ref{sequence-lambda})}, then for any
 positive integer $k$ and real vector
 $(v_1,v_2,\ldots,v_k)^T\in\mathbb{R}^k$,
 it holds that
 $$\dsum_{n=0}^{k-1}\Big(\dsum_{p=0}^n\lambda_p^{(\alpha)}v_{n+1-p}\Big)v_{n+1}\geq0.$$
 \end{lemma}
 {\bf Proof. } For simplicity of presentation, in this proof,
 we denote $g_p=g_p^{(\alpha)},~\lambda_p=\lambda_p^{(\alpha)}$
 without ambiguity. One can easily check that,
 to prove the above quadratic form is nonnegative is equivalent to proving
 the symmetric Toeplitz matrix T is positive semi-definite, where
 $$T=\left(
 \renewcommand{\arraystretch}{0.8}
 \begin{array}{ccccc}
 \lambda_0 & \frac{\lambda_1}2 & \frac{\lambda_2}2 & \cdots & \frac{\lambda_{k-1}}2 \\
 \frac{\lambda_1}2 & \lambda_0 & \frac{\lambda_1}2 & \ddots & \vdots \\
 \frac{\lambda_2}2 & \frac{\lambda_1}2 & \ddots & \ddots & \frac{\lambda_2}2 \\
 \vdots & \ddots & \ddots & \lambda_0 & \frac{\lambda_1}2 \\
 \frac{\lambda_{k-1}}2 & \cdots & \frac{\lambda_2}2 & \frac{\lambda_1}2 & \lambda_0
\end{array}
\right).$$
 Notice that the generating function (see \cite{Chan2}) of T is given by
 $$\begin{array}{rl}
 f(\alpha,x)&=\lambda_0+\dfrac12\dsum_{k=1}^\infty\lambda_ke^{ikx}
 +\dfrac12\dsum_{k=1}^\infty\lambda_ke^{-ikx}\\
 &=\dfrac{2+\alpha}2{g}_0
 +\dfrac12\dsum_{k=1}^\infty\bigg(\dfrac{2+\alpha}2{g}_k-\dfrac{\alpha}2{g}_{k-1}\bigg)e^{ikx}
 +\dfrac12\dsum_{k=1}^\infty\bigg(\dfrac{2+\alpha}2{g}_k-\dfrac{\alpha}2{g}_{k-1}\bigg)e^{-ikx}\\
 &=\dfrac{2+\alpha}4(1-e^{ix})^{\alpha}
 +\dfrac{2+\alpha}4(1-e^{-ix})^{\alpha}
 -\dfrac{\alpha}4e^{ix}(1-e^{ix})^{\alpha}
 -\dfrac{\alpha}4e^{-ix}(1-e^{-ix})^{\alpha}.
 \end{array}$$
 As mentioned in \cite{Tian}, we  only need to consider the principal value of $f(\alpha,x)$
 on $[0, \pi]$ which gives
 $$\begin{array}{rl}
 f(\alpha,x)&=\dfrac{2+\alpha}4[2i\sin(-\dfrac{x}2)e^{\frac{ix}2}]^{\alpha}
 +\dfrac{2+\alpha}4[2i\sin(\dfrac{x}2)e^{\frac{-ix}2}]^{\alpha}\\[7pt]
 &\quad-\dfrac{\alpha}4e^{ix}[2i\sin(-\dfrac{x}2)e^{\frac{ix}2}]^{\alpha}
 -\dfrac{\alpha}4e^{-ix}[2i\sin(\dfrac{x}2)e^{\frac{-ix}2}]^{\alpha}\\[7pt]
 &=[2\sin(\dfrac{x}2)]^{\alpha}\Big\{\dfrac{2+\alpha}2\cos[\dfrac{\alpha}2(\pi-x)]
 -\dfrac{\alpha}2\cos[\dfrac{\alpha}2(\pi-x)-x]\Big\}.
 \end{array}$$
 Let $h(\alpha,x)=\dfrac{2+\alpha}2\cos[\dfrac{\alpha}2(\pi-x)]
 -\dfrac{\alpha}2\cos[\dfrac{\alpha}2(\pi-x)-x]$, then one can
 easily check that
 $$h_x(\alpha,x)=\dfrac{\alpha(2+\alpha)}2
 \sin(\dfrac{x}2)\cos[\dfrac{\alpha}2(\pi-x)-\frac x2]
 \geq0.$$
 Therefore $h(\alpha,x)$ is nondecreasing with respect to $x$ and
 $h(\alpha,x)\geq h(\alpha,0)=\cos(\frac{\alpha\pi}2)\geq0$, which implies that
 $f(\alpha,x)\geq0$.
 The lemma now follows as a result of the Grenander-Szeg\"{o} Theorem \cite{Chan2}. $\qquad\Box$

 \begin{remark}
 We note here that the function $h(\alpha,x)$ can not be proved to be nonnegative by considering differentiation with respect to $\alpha$ as in
 {\rm\cite{Tian}}.
 \end{remark}

 \begin{lemma}\label{gronwall}{\rm(Grownall's inequality \cite{Quarteroni})}
 Assume that $\{k_n\}$ and $\{p_n\}$ are nonnegative sequences, and the sequence $\{\phi_n\}$ satisfies
 $$\phi_0\leq g_0,~~~\phi_n\leq g_0+\dsum_{l=0}^{n-1}p_l
 +\dsum_{l=0}^{n-1}k_l\phi_l,~~~n\geq1,$$
 where $g_0\geq0$. Then the sequence $\{\phi_n\}$ satisfies
 $$\phi_n\leq\Big(g_0+\dsum_{l=0}^{n-1}p_l\Big)
 \exp{\Big(\dsum_{l=0}^{n-1}k_l\Big)},~~~n\geq1.$$
 \end{lemma}

 With all the preparation,
 we can now show the convergence and stability of the compact finite difference scheme (\ref{compact-scheme1})--(\ref{compact-scheme2}).

 \begin{theorem}\label{sub-convergence}
 Assume that $u(x,t)\in {\cal C}_{x,t}^{6,2}([0,L]\times[0,T])$ is the solution of {\rm(\ref{sub-main1})--(\ref{sub-main2})}
  and $\{u_j^k|0\leq j\leq M,~ 0\leq k\leq N\}$
 is the solution of the finite difference scheme {\rm(\ref{compact-scheme1})--(\ref{compact-scheme2})}, respectively.
 Denote
 $$e_i^k=u(x_i,t_k)-u_i^k, \quad 0\leq i\leq M, \quad 0\leq k\leq N.$$
 Then there exists a positive constant $\tilde c_1$ such that
 $$\|e^k\|\leq\tilde c_1(\tau^2+h^4), \quad 0\leq k\leq N.$$
 \end{theorem}
 {\bf Proof.} We can easily get the error equation:
 \begin{equation}\label{error}
 \begin{array}{l}
 {\cal H}(e_i^{k+1}-e_i^{k})
 =\dfrac{\kappa_1\tau^{1-\alpha}}2
 \dsum_{l=0}^{k}\lambda_l^{(\alpha)}\delta^2_x(e^{k+1-l}_i+e^{k-l}_i)
 +\dfrac{\kappa_2\tau^{1-\beta}}2
 \dsum_{l=0}^{k}\lambda_l^{(\beta)}\delta^2_x(e^{k+1-l}_i+e^{k-l}_i)
 +\tau R_i^{k+1}\\
 0\leq k\leq N-1,~~1\leq i\leq M-1,\\
 e_0^k=e_M^k=0, \quad 1\leq k\leq N\\
 e_i^0=0, \quad 0\leq i\leq M,
 \end{array}
 \end{equation}
 where $|R_i^{k+1}|\leq c_1(\tau^{2}+h^4)$.

 Denote $C=\frac1{12}tri[1,10,1]$, multiplying (\ref{error}) by $h(e_i^{k+1}+e_i^{k})$ and summing in $i$, we obtain
 \begin{equation}\label{for-stability}
 \begin{array}{rl}
 h(e^{k+1}+e^{k})^TC(e^{k+1}-e^{k})
 =&\dfrac{\kappa_1\tau^{1-\alpha}}2\dsum_{l=0}^{k}\lambda_l^{(\alpha)}\langle\delta_x^2(e^{k+1-l}+e^{k-l}),e^{k+1}+e^{k}\rangle\\
 &+\dfrac{\kappa_2\tau^{1-\beta}}2\dsum_{l=0}^{k}\lambda_l^{(\beta)}\langle\delta_x^2(e^{k+1-l}+e^{k-l}),e^{k+1}+e^{k}\rangle\\
 & +\tau h(e^{k+1}+e^{k})^TR^{k+1}\\
 =&-\dfrac{\kappa_1\tau^{1-\alpha}}2\dsum_{l=0}^{k}\lambda_l^{(\alpha)}\langle\delta_x(e^{k+1-l}+e^{k-l}),\delta_x(e^{k+1}+e^{k})\rangle\\
 &-\dfrac{\kappa_2\tau^{1-\beta}}2\dsum_{l=0}^{k}\lambda_l^{(\beta)}\langle\delta_x(e^{k+1-l}+e^{k-l}),\delta_x(e^{k+1}+e^{k})\rangle\\
 &+\tau h(e^{k+1}+e^{k})^TR^{k+1}.
 \end{array}
 \end{equation}
 Summing up for $0\leq k\leq n-1$ and noticing that
 $$h(e^{k+1}+e^{k})^TC(e^{k+1}-e^{k})=h({e^{k+1}}^TCe^{k+1}-{e^k}^TCe^{k}),$$
 $$h{e^n}^TCe^n\geq\dfrac23\|e^n\|^2,$$
 and
 $$\tau h(e^{k+1}+e^{k})^TR^{k+1} \leq\dfrac\tau3(\|e^{k+1}\|^2+\|e^{k}\|^2)+\dfrac{3\tau}2\|R^{k+1}\|^2,$$
 we get, by Lemma \ref{sub-sequence}, that
 $$\begin{array}{rl}
 \dfrac23\|e^n\|^2
 \leq&-\dfrac{\kappa_1\tau^{1-\alpha}}2\dsum_{k=0}^{n-1}\dsum_{l=0}^{k}\lambda_l^{(\alpha)}\langle\delta_x(e^{k+1-l}+e^{k-l}),\delta_x(e^{k+1}+e^{k})\rangle\\
 &-\dfrac{\kappa_2\tau^{1-\beta}}2\dsum_{k=0}^{n-1}\dsum_{l=0}^{k}\lambda_l^{(\beta)}\langle\delta_x(e^{k+1-l}+e^{k-l}),\delta_x(e^{k+1}+e^{k})\rangle\\
 &+\dfrac\tau3\dsum_{k=0}^{n-2}(\|e^{k+1}\|^2+\|e^{k}\|^2)
 +\dfrac{3\tau}2\dsum_{k=0}^{n-2}\|R^{k+1}\|^2
 +\tau h(e^{n}+e^{n-1})^TR^{n}\\
 \leq&\dfrac13\|e^n\|^2+\dfrac\tau3\|e^{n-1}\|^2
 +\dfrac{3\tau^2}4\|R^{n}\|^2+\dfrac{3\tau}4\|R^{n}\|^2\\[5pt]
 &+\dfrac\tau3\dsum_{k=1}^{n-1}\|e^{k}\|^2+\dfrac\tau3\dsum_{k=1}^{n-2}\|e^{k}\|^2
 +\dfrac{3\tau}2\dsum_{k=0}^{n-2}\|R^{k+1}\|^2\\
 \leq&\dfrac13\|e^n\|^2+\dfrac{3\tau^2}4\|R^{n}\|^2
 +\dfrac{2\tau}3\dsum_{k=1}^{n-1}\|e^{k}\|^2
 +\dfrac{3\tau}2\dsum_{k=0}^{n-1}\|R^{k+1}\|^2,
 \end{array}$$
 which gives
 $$\begin{array}{rl}
 \|e^n\|^2&\leq\dfrac{9\tau^2}4\|R^{n}\|^2
 +2\tau\dsum_{k=1}^{n-1}\|e^{k}\|^2
 +\dfrac{9\tau}2\dsum_{k=0}^{n-1}\|R^{k+1}\|^2\\
 &\leq2\tau\dsum_{k=1}^{n-1}\|e^{k}\|^2+c(\tau^2+h^4)^2.
 \end{array}$$
 then the desired result follows by Lemma \ref{gronwall}.$\qquad\Box$

 \begin{remark}\label{sub-stability}
 By using similar techniques, we can show that the compact scheme {\rm(\ref{compact-scheme1})--(\ref{compact-scheme2})} is stable for $0<\alpha,~\beta<1$.
 In fact, suppose that $\{v_i^l\}$ is the solution of
 \begin{align}\label{pur1}\nonumber
 &{\cal H}(v^{k+1}_i-v^k_i)=\dfrac{\kappa_1\tau^{1-\alpha}}2
 \bigg(\dsum_{l=0}^{k+1}\lambda_l^{(\alpha)}\delta^2_xv^{k+1-l}_i
 +\dsum_{l=0}^{k}\lambda_l^{(\alpha)}\delta^2_xv^{k-l}_i\bigg)\\\nonumber
 &\quad+\dfrac{\kappa_2\tau^{1-\beta}}2
 \bigg(\dsum_{l=0}^{k+1}\lambda_l^{(\beta)}\delta^2_xv^{k+1-l}_i
 +\dsum_{l=0}^{k}\lambda_l^{(\beta)}\delta^2_xv^{k-l}_i\bigg)
 +\dfrac\tau2{\cal H}(f^k_i+f^{k+1}_i),\\
 &0\leq k\leq N-1,~~1\leq i\leq M-1,\\\nonumber
 &v_0^k=\varphi_1^k,\quad v_M^k=\varphi_2^k,\quad 1\leq k\leq N,\\\label{pur2}
 &v_i^0=\rho_i,\quad 0\leq i\leq M,
 \end{align}
 where $\rho\in {\cal V}$.
 Then, by subtracting {\rm(\ref{compact-scheme1})--(\ref{compact-scheme2})}
 from {\rm(\ref{pur1})--(\ref{pur2})}, we obtain the following equations for  $\varepsilon_i^l=v_i^l-u_i^l-\rho_i$,
 \begin{align}\label{pur3}\nonumber
 &{\cal H}(\varepsilon_i^{k+1}-\varepsilon_i^{k})
 =\dfrac{\kappa_1\tau^{1-\alpha}}2
 \bigg(\dsum_{l=0}^{k+1}\lambda_l^{(\alpha)}\delta^2_x\varepsilon^{k+1-l}_i
 +\dsum_{l=0}^{k}\lambda_l^{(\alpha)}\delta^2_x\varepsilon^{k-l}_i\bigg)
 +\frac{\kappa_1\tau^{1-\alpha}}2\bigg(\dsum_{l=0}^{k+1}\lambda_l^{(\alpha)}
 +\dsum_{l=0}^{k}\lambda_l^{(\alpha)}\bigg)\delta^2_x\rho_i\\
 &\quad+\dfrac{\kappa_2\tau^{1-\beta}}2
 \bigg(\dsum_{l=0}^{k+1}\lambda_l^{(\beta)}\delta^2_x\varepsilon^{k+1-l}_i
 +\dsum_{l=0}^{k}\lambda_l^{(\beta)}\delta^2_x\varepsilon^{k-l}_i\bigg)
 +\frac{\kappa_2\tau^{1-\beta}}2\bigg(\dsum_{l=0}^{k+1}\lambda_l^{(\beta)}
 +\dsum_{l=0}^{k}\lambda_l^{(\beta)}\bigg)\delta^2_x\rho_i,\\\nonumber
 &0\leq k\leq N-1,~~1\leq i\leq M-1,\\\nonumber
 &\varepsilon_0^k=\varepsilon_M^k=0,\quad 1\leq k\leq N,\\\nonumber
 &\varepsilon_i^0=0,\quad 0\leq i\leq M.
 \end{align}
 Notice that {\rm\cite{r2,r3}}, $g_l^{(\alpha)}$ and $g_l^{(\beta)}$ are less than 0 for $l\geq1$ and $$\tau^{-\alpha}\sum\limits_{l=0}^{k}g_l^{(\alpha)}=\frac1{\Gamma(1-\alpha)}+O(\tau),
 ~\tau^{-\beta}\sum\limits_{l=0}^{k}g_l^{(\beta)}=\frac1{\Gamma(1-\beta)}+O(\tau).$$
 We can therefore obtain
 \begin{align*}
 \frac{\tau^{-\alpha}}2\bigg(\dsum_{l=0}^{k+1}\lambda_l^{(\alpha)}
 +\dsum_{l=0}^{k}\lambda_l^{(\alpha)}\bigg)
 &=\frac{\tau^{-\alpha}}2\bigg[\Big(1+\frac{\alpha}2\Big)\sum\limits_{l=0}^{k+1}g_l^{(\alpha)}
 +\sum\limits_{l=0}^{k}g_l^{(\alpha)}-\frac{\alpha}2\sum\limits_{l=0}^{k-1}g_{l}^{(\alpha)}\bigg]\\
 &=\frac1{\Gamma(1-\alpha)}+O(\tau),
 \end{align*}
 and
 $$\frac{\tau^{-\beta}}2\bigg(\dsum_{l=0}^{k+1}\lambda_l^{(\beta)}
 +\dsum_{l=0}^{k}\lambda_l^{(\beta)}\bigg)
 =\frac1{\Gamma(1-\beta)}+O(\tau).$$
 Now multiplying {\rm(\ref{pur3})} by $h(\varepsilon_i^{k+1}+\varepsilon_i^{k})$
 and summing in $i$ and $k$, by  arguments similar to that for the proof of Theorem {\rm\ref{sub-convergence}},
 we can get that
 $$\|\varepsilon^k\|^2\leq2\tau\dsum_{l=0}^{k-1}\|\varepsilon^{l}\|^2+5c_2^2\|\delta_x^2\rho\|^2,$$
 where $c_2=1+\frac{\kappa_1}{\Gamma(1-\alpha)}+\frac{\kappa_2}{\Gamma(1-\beta)}$,
 it then follows from Lemma {\rm\ref{gronwall}} that $\|\varepsilon^k\|^2\leq5c_2^2e^{2T}\|\delta_x^2\rho\|^2$.
 Finally we can conclude that $\|v^k-u^k\|\leq\|v^k-u^k-\rho\|+\|\rho\|\leq\sqrt5c_2e^{T}\|\delta_x^2\rho\|+\|\rho\|$.
 \end{remark}
 \medskip
 Before the numerical experiments, we first turn to the study on the fractional diffusion-wave equation.

 \section{A Compact scheme for the fractional diffusion-wave equation}
 In this section, we consider the following time fractional diffusion-wave equation
 \begin{equation}\label{dw-main1}
 \begin{array}{rl}
 &_{0}^CD_t^{\gamma}u
 =\kappa\dfrac{\partial^2u}{\partial x^2}+g(x,t),
  \quad 0\leq x\leq L, \quad 0<t\leq T, \quad 1<\gamma<2,\\
  \end{array}
 \end{equation}
 subject to
 \begin{equation}\label{dw-main2}
 \begin{array}{rl}
 &u(x,0)=0, \quad \dfrac{\partial u(x,0)}{\partial t}=\phi(x), \quad 0\leq x\leq L,\\
 &u(0,t)=\varphi_1(t), \quad u(L,t)=\varphi_2(t), \quad 0<t\leq T,
 \end{array}
 \end{equation}
 where $\kappa$ is a positive constant, $_{0}^CD_t^{\gamma}u$ is the Caputo
 fractional derivative of $u$ with respect to time variable $t$.

 However, instead of solving (\ref{dw-main1}) directly,
 we follow the technique in \cite{Huang},
 where the problem is equivalently changed to the following:
 \begin{equation}\label{changed-main1}
 \dfrac{\partial u(x,t)}{\partial t}=\phi(x)+
 \dfrac\kappa{\Gamma(\alpha)}\dint_0^t(t-s)^{\alpha-1}\dfrac{\partial^2u(x,s)}{\partial x^2}ds+f(x,t), \quad 0\leq x\leq L, \quad 0<t\leq T,
 \end{equation}
 where $0<\alpha=\gamma-1<1,~f(x,t)={_{0}I^\alpha_tg(x,t)}$, and $_{0}I^\alpha_t$ is the Riemann-Liouville fractional integral operator
 with respect to $t$.

 Once again, we choose $(p,q)=(0,-1)$ in (ii) of Lemma \ref{main-operator},
 yielding ${\mu_1}=1-\frac\alpha2,~{\mu_2}=\frac\alpha2$ in (\ref{weight}),  and
 $$\begin{array}{rl}
 _{0}I^\alpha_tu_{xx}(x_i,t_{n+1})
 &=\tau^\alpha\bigg[\Big(1-\dfrac\alpha2\Big)\dsum_{k=0}^{n+1}\omega_k^{(\alpha)}\delta^2_xu^{n+1-k}_i
 +\dfrac\alpha2\dsum_{k=0}^{n}\omega_k^{(\alpha)}\delta^2_xu^{n-k}_i\bigg]+O(\tau^2+h^2)\\
 &=\tau^\alpha\dsum_{k=0}^{n+1}\lambda_k\delta^2_xu^{n+1-k}_i+O(\tau^2+h^2),
 \end{array}$$
 where
 \begin{equation}\label{wave-sequence}
 \lambda_0=(1-\dfrac\alpha2)\omega_0^{(\alpha)},
 ~~\lambda_k=(1-\dfrac\alpha2)\omega_{k}^{(\alpha)}
 +\dfrac\alpha2\omega_{k-1}^{(\alpha)},~~k\geq1.
 \end{equation}
 We can therefore introduce the compact scheme of (\ref{dw-main2})--(\ref{changed-main1}) as
 \begin{align}\label{dw-compact-scheme1}\nonumber
 &{\cal H}(u^{n+1}_i-u^n_i)=\tau{\cal H}\phi_i+\dfrac{\kappa\tau^{\alpha+1}}2
 \bigg(\dsum_{k=0}^{n+1}\lambda_k\delta^2_xu^{n+1-k}_i
 +\dsum_{k=0}^{n}\lambda_k\delta^2_xu^{n-k}_i\bigg)
 +\dfrac\tau2{\cal H}(f^n_i+f^{n+1}_i),\\
 &~0\leq n\leq N-1,~~1\leq i\leq M-1,\\\nonumber
 &u_0^n=\varphi_1^n, \quad u_M^n=\varphi_2^n, \quad 1\leq n\leq N,\\\label{dw-compact-scheme2}
 &u_i^0=0,\quad 0\leq i\leq M.
 \end{align}

 Due to the same reasons stated in the last section, the difference scheme (\ref{dw-compact-scheme1})--(\ref{dw-compact-scheme2}) has a unique solution.

 Similar to the role of Lemma \ref{sub-sequence} to Theorem \ref{sub-convergence}, convergence of the  scheme (\ref{dw-compact-scheme1})--(\ref{dw-compact-scheme2})
 is established using the following lemma:
 \begin{lemma}\label{sequence}
 Let $\{\lambda_n\}_{n=0}^\infty$ be defined as in {\rm(\ref{wave-sequence})}, then for any
 positive integer $k$ and real vector
 $(v_1,v_2,\ldots,v_k)^T\in\mathbb{R}^k$,
 it holds that
 $$\dsum_{n=0}^{k-1}\Big(\dsum_{p=0}^n\lambda_pv_{n+1-p}\Big)v_{n+1}\geq0.$$
 \end{lemma}
 {\bf Proof. }
 The generating function of the Toeplitz matrix for the corresponding sequence is given by
 $$\begin{array}{rl}
 f(\alpha,x)&=\lambda_0+\dfrac12\dsum_{k=1}^\infty\lambda_ke^{ikx}
 +\dfrac12\dsum_{k=1}^\infty\lambda_ke^{-ikx}\\
 &=(1-\dfrac\alpha2)\omega_0^{(\alpha)}
 +\dfrac12\dsum_{k=1}^\infty\big[(1-\dfrac\alpha2)\omega_k^{(\alpha)}
 +\dfrac\alpha2\omega_{k-1}^{(\alpha)}\big]e^{ikx}
 +\dfrac12\dsum_{k=1}^\infty\big[(1-\dfrac\alpha2)\omega_k^{(\alpha)}
 +\dfrac\alpha2\omega_{k-1}^{(\alpha)}\big]e^{-ikx}\\[7pt]
 &=\dfrac{2-\alpha}4(1-e^{ix})^{-\alpha}
 +\dfrac{2-\alpha}4(1-e^{-ix})^{-\alpha}
 +\dfrac{\alpha}4e^{ix}(1-e^{ix})^{-\alpha}
 +\dfrac{\alpha}4e^{-ix}(1-e^{-ix})^{-\alpha}\\[7pt]
 &=[2\sin(\dfrac{x}2)]^{-\alpha}\Big\{(1-\dfrac\alpha2)\cos[\dfrac\alpha2(\pi-x)]
 +\dfrac\alpha2\cos[x+\dfrac\alpha2(\pi-x)]\Big\}.
 \end{array}$$
 Let $g(\alpha,x)=(1-\dfrac\alpha2)\cos[\dfrac\alpha2(\pi-x)]
 +\dfrac\alpha2\cos[x+\dfrac\alpha2(\pi-x)]$, then we can easily compute
 $$g_\alpha=-\sin(\dfrac{x}2)\sin[\dfrac{x}2+\dfrac\alpha2(\pi-x)]
 -(1-\dfrac\alpha2)\sin[\dfrac\alpha2(\pi-x)]\dfrac{\pi-x}2
 -\dfrac\alpha2\sin[x+\dfrac\alpha2(\pi-x)]\dfrac{\pi-x}2
 \leq0.$$
 Since $g(\alpha,x)$ is nonincreasing with respect to $\alpha$,
 we have $g(\alpha,x)\geq g(1,x)=0$. This yields $f(\alpha,x)\geq0$ and the result follows. $\qquad\Box$

 \begin{remark}
 {\rm(i)} We point out  that interpolating polynomials have been used
 in {\rm\cite{Li}} to approximate the Riemann-Liouville fractional
 integral with second order accuracy. However, one can easily test
 that the corresponding coefficients do not satisfy Lemma
 {\rm\ref{sequence}} in general.

 {\rm(ii)} The Grenander-Szeg\"{o} Theorem is proved for continuous
 generating function in {\rm\cite{Chan2}}. However, one can check
 that the arguments used in the proof of the  Grenander-Szeg\"{o} Theorem can still be applied here
 to conclude Lemma {\rm\ref{sequence}}.
 \end{remark}

 With Lemma \ref{sequence}, we can obtain the following convergence result of our compact difference scheme (\ref{dw-compact-scheme1})--(\ref{dw-compact-scheme2}).
 Since the proof is similar to that of Theorem \ref{sub-convergence}, we therefore skip the details.
 \begin{theorem}
 Assume that $u(x,t)\in {\cal C}_{x,t}^{6,2}([0,L]\times[0,T])$ is the solution of
 {\rm(\ref{dw-main2})}--{\rm(\ref{changed-main1})}
 and $\{u_i^k|0\leq i\leq M,~ 0\leq k\leq N\}$
 is the solution of the finite difference scheme {\rm(\ref{dw-compact-scheme1})}--{\rm(\ref{dw-compact-scheme2})}, respectively. Denote $$e_i^k=u(x_i,t_k)-u_i^k, \quad 0\leq i\leq M, \quad 0\leq k\leq N.$$
 Then there exists a positive constant $\tilde c_2$ such that
 $$\|e^k\|\leq\tilde c_2(\tau^2+h^4), \quad 0\leq k\leq N.$$
 \end{theorem}
 \begin{remark}
 Following the ideas of the proof for Theorem {\rm\ref{sub-convergence}}
 and Remark {\rm\ref{sub-stability}},
 one can show that the proposed compact scheme {\rm(\ref{dw-compact-scheme1})--(\ref{dw-compact-scheme2})} is stable
  in the sense that if $\{v_i^l\}$ is the solution of
 \begin{align*}
 &{\cal H}(v^{k+1}_i-v^k_i)=\tau{\cal H}(\phi_i+\tilde\rho_i)+\dfrac{\kappa\tau^{\alpha+1}}2
 \bigg(\dsum_{l=0}^{k+1}\lambda_l\delta^2_xv^{k+1-l}_i
 +\dsum_{l=0}^{k}\lambda_l\delta^2_xv^{k-l}_i\bigg)
 +\dfrac\tau2{\cal H}(f^k_i+f^{k+1}_i),\\
 &~0\leq k\leq N-1,~~1\leq i\leq M-1,\\
 &v_0^k=\varphi_1^k,\quad v_M^k=\varphi_2^k,\quad 1\leq k\leq N,\\
 &v_i^0=\rho_i,\quad 0\leq i\leq M,
 \end{align*}
 with $\rho, \tilde\rho\in {\cal V}$,
 then
 $$\|v^k-u^k\|\leq\sqrt5e^{T}\Big\|\Big(\frac{\kappa}{\Gamma(\alpha+1)}+1\Big)\delta_x^2\rho+\tilde\rho\Big\|+\|\rho\|.$$
 \end{remark}

 \section{Numerical experiments}
 In this section, we carry out numerical experiments using the proposed finite difference schemes
 to illustrate our theoretical statements. All our tests were done in MATLAB. We suppose $L=T=1$.
 The maximum norm errors and 2-norm errors between the exact and the numerical solutions
 $$E_\infty(h,\tau)=\max_{0\leq k\leq N}\|U^k-u^k\|_\infty, \quad E_2(h,\tau)=\max_{0\leq k\leq N}\|U^k-u^k\|$$
 are shown.
 Furthermore, the temporal convergence order, denoted by
 $$Rate1_\infty=\log_2\bigg(\dfrac{E_\infty(h,2\tau)}{E_\infty(h,\tau)}\bigg),
 \quad Rate1_2=\log_2\bigg(\dfrac{E_2(h,2\tau)}{E_2(h,\tau)}\bigg)$$
 for sufficiently small $h$, and the spatial convergence order, denoted
 by
 $$Rate2_\infty=\log_2\bigg(\dfrac{E_\infty(2h,\tau)}{E_\infty(h,\tau)}\bigg),
 \quad Rate2_2=\log_2\bigg(\dfrac{E_2(2h,\tau)}{E_2(h,\tau)}\bigg)$$
 when $\tau$ is sufficiently small, are reported.
 The numerical results given by these examples justify our theoretical
 analysis.

 \bigskip
 \begin{example}\label{ex2}
 The following problem was studied in
 {\rm\cite{Mohebbi}}:
 $$\begin{array}{rl}
 &\dfrac{\partial u(x,t)}{\partial t}
 =\bigg(\dfrac{\partial^{1-\tilde\alpha}}{\partial t^{1-\tilde\alpha}}
 +\dfrac{\partial^{1-\tilde\beta}}{\partial t^{1-\tilde\beta}}\bigg)
 \bigg[\dfrac{\partial^2u(x,t)}{\partial x^2}\bigg]+g(x,t),\\
 &u(0,t)=0,~~u(1,t)=t^{1+\tilde\alpha+\tilde\beta}\sin(1),~~0\leq t\leq 1,\\
 &u(x,0)=0,~~0\leq x\leq 1,
 \end{array}$$
 where
 $$g(x,t)=\sin(x)\bigg[(1+\tilde\alpha+\tilde\beta)t^{\tilde\alpha+\tilde\beta}
 +\frac{\Gamma(2+\tilde\alpha+\tilde\beta)}{\Gamma(1+2\tilde\alpha+\tilde\beta)}t^{2\tilde\alpha+\tilde\beta}
 +\frac{\Gamma(2+\tilde\alpha+\tilde\beta)}{\Gamma(1+\tilde\alpha+2\tilde\beta)}t^{\tilde\alpha+2\tilde\beta}\bigg].$$
 In terms of  the notations in this paper, we consider
 $$\begin{array}{rl}
 &\dfrac{\partial u(x,t)}{\partial t}
 =\bigg(\dfrac{\partial^{\alpha}}{\partial t^{\alpha}}
 +\dfrac{\partial^{\beta}}{\partial t^{\beta}}\bigg)
 \bigg[\dfrac{\partial^2u(x,t)}{\partial x^2}\bigg]+f(x,t),\\
 &u(0,t)=0,~~u(1,t)=t^{3-\alpha-\beta}\sin(1),~~0\leq t\leq 1,\\
 &u(x,0)=0,~~0\leq x\leq 1,
 \end{array}$$
 with
 $$f(x,t)=\sin(x)\bigg[(3-\alpha-\beta)t^{2-\alpha-\beta} +\frac{\Gamma(4-\alpha-\beta)}{\Gamma(4-2\alpha-\beta)}t^{3-2\alpha-\beta}
 +\frac{\Gamma(4-\alpha-\beta)}{\Gamma(4-\alpha-2\beta)}t^{3-\alpha-2\beta}\bigg].$$
 It is easy to check that the exact solution is $u(x,t)=\sin(x)t^{3-\alpha-\beta}$.
 \end{example}

 \begin{figure}\label{figure1}
 \begin{center}
 \includegraphics[height=6.3cm,width=17.3cm]{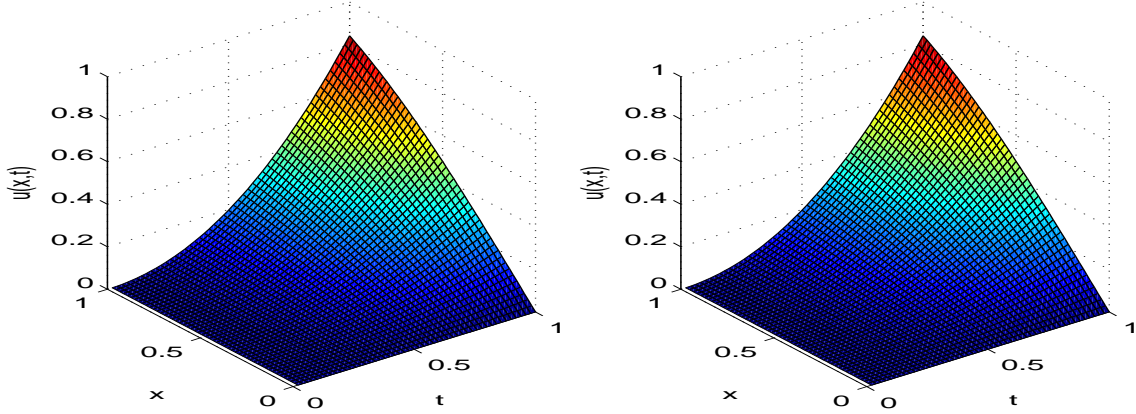}
 \end{center}
 \caption{The exact solution (left) and numerical solution (right) for Example \ref{ex2}, when $\alpha=0.2,~\beta=0.7,~h=\tau=\frac1{50}$.}
 \end{figure}

 \begin{table}[h]
 \begin{center}
 \caption{Numerical convergence orders in temporal direction with $h=\frac1{30},~\alpha=0.35,~\beta=0.05$ for Example \ref{ex2}.}

 \bigskip
 \label{table3}
 \noindent
 \renewcommand{\arraystretch}{0.9}
 {\small
 \begin{tabular}{lllll}
 \hline
 $\tau$   & $E_\infty(h,\tau)$ & $Rate1_\infty$& $E_2(h,\tau)$ & $Rate1_2$\\\hline
 1/5      & 8.2083e-4          & $\ast$        & 5.9004e-4     & $\ast$\\
 1/10     & 1.9894e-4          & 2.0447        & 1.4192e-4     & 2.0557\\
 1/20     & 5.0072e-5          & 1.9903        & 3.5877e-5     & 1.9840\\
 1/40     & 1.2567e-5          & 1.9944        & 9.0135e-6     & 1.9929\\
 1/80     & 3.1502e-6          & 1.9961        & 2.2593e-6     & 1.9962\\
 1/160    & 7.8883e-7          & 1.9976        & 5.6574e-7     & 1.9976\\
 \hline
 \end{tabular}
 }
 \end{center}
 \end{table}

 \begin{table}[h]
 \begin{center}
 \caption{Numerical convergence orders in temporal direction with $h=\frac1{30},~\alpha=0.2,~\beta=0.7$ for Example \ref{ex2}.}

 \bigskip
 \label{table4}
 \noindent
 \renewcommand{\arraystretch}{0.9}
 {\small
 \begin{tabular}{lllll}
 \hline
 $\tau$   & $E_\infty(h,\tau)$ & $Rate1_\infty$& $E_2(h,\tau)$ & $Rate1_2$\\\hline
 1/5      & 7.3332e-4          & $\ast$        & 5.2626e-4     & $\ast$\\
 1/10     & 1.8938e-4          & 1.9532        & 1.3595e-4     & 1.9527\\
 1/20     & 4.8407e-5          & 1.9680        & 3.4765e-5     & 1.9674\\
 1/40     & 1.2305e-5          & 1.9759        & 8.8421e-6     & 1.9752\\
 1/80     & 3.1336e-6          & 1.9734        & 2.2500e-6     & 1.9745\\
 1/160    & 8.0171e-7          & 1.9667        & 5.7465e-7     & 1.9692\\
 \hline
 \end{tabular}
 }
 \end{center}
 \end{table}

 \begin{table}[h]
 \begin{center}
 \caption{Numerical convergence orders in spatial direction with $\tau=\frac1{8000},~\alpha=0.15,~\beta=0.35$ for Example \ref{ex2}.}

 \bigskip
 \label{table5}
 \noindent
 \renewcommand{\arraystretch}{0.9}
 {\small
 \begin{tabular}{lllll}
 \hline
 $h$   & $E_\infty(h,\tau)$ & $Rate2_\infty$& $E_2(h,\tau)$ & $Rate2_2$\\\hline
 1/2   & 1.1910e-5          & $\ast$        & 9.8357e-6     & $\ast$\\
 1/4   & 8.6322e-7          & 4.0102        & 6.3362e-7     & 3.9563\\
 1/8   & 5.4954e-8          & 3.9734        & 3.9832e-8     & 3.9916\\
 1/16  & 3.7617e-9          & 3.8688        & 2.7001e-9     & 3.8829\\
 \hline
 \end{tabular}
 }
 \end{center}
 \end{table}

 \begin{table}[h]
 \begin{center}
 \caption{Numerical convergence orders in spatial direction with $\tau=\frac1{8000},~\alpha=0.75,~\beta=0.15$ for Example \ref{ex2}.}

 \bigskip
 \label{table6}
 \noindent
 \renewcommand{\arraystretch}{0.9}
 {\small
 \begin{tabular}{lllll}
 \hline
 $h$   & $E_\infty(h,\tau)$ & $Rate2_\infty$& $E_2(h,\tau)$ & $Rate2_2$\\\hline
 1/2   & 1.4250e-5          & $\ast$        & 1.0077e-5     & $\ast$\\
 1/4   & 8.8434e-7          & 4.0103        & 6.4856e-7     & 3.9576\\
 1/8   & 5.6182e-8          & 3.9764        & 4.0749e-8     & 3.9924\\
 1/16  & 3.8214e-9          & 3.8779        & 2.7425e-9     & 3.8932\\
 \hline
 \end{tabular}
 }
 \end{center}
 \end{table}

 Figure 1 shows the exact solution (left) and numerical solution (right) for Example \ref{ex2},
 when $\alpha=0.2,~\beta=0.7,~h=\tau=\frac1{50}$. Meanwhile,
 we list, in Table \ref{table3} and Table \ref{table4}, the convergence order in temporal direction with $h=\frac1{30}$
 and, in Table \ref{table5} and Table \ref{table6}, the convergence order in spatial direction with $\tau=\frac1{8000}$.
 The convergence order of the numerical results matches that of the theoretical one.

 \bigskip
 We next consider a special case of the modified anomalous fractional sub-diffusion
 equation:
 \begin{example}\label{ex1}
 Consider the following  example {\rm(\cite{Cui,Lin})}.
 $$\begin{array}{rl}
 &\dfrac{\partial u(x,t)}{\partial t}
 =\dfrac{\partial^{1-\tilde\alpha}}{\partial t^{1-\tilde\alpha}}
 \bigg[\dfrac{\partial^2u(x,t)}{\partial x^2}\bigg]
 +\bigg[2t+\frac{8\pi^2t^{1+\tilde\alpha}}{\Gamma(2+\tilde\alpha)}\bigg]\sin(2\pi x),\\
 &u(0,t)=u(1,t)=0,~~0\leq t\leq 1,\\
 &u(x,0)=0,~~0\leq x\leq 1.
 \end{array}$$
  The problem can be written as
 $$\begin{array}{rl}
 &\dfrac{\partial u(x,t)}{\partial t}
 =\dfrac{\partial^{\alpha}}{\partial t^{\alpha}}
 \bigg[\dfrac{\partial^2u(x,t)}{\partial x^2}\bigg]
 +\bigg[2t+\frac{8\pi^2t^{2-\alpha}}{\Gamma(3-\alpha)}\bigg]\sin(2\pi x),\\
 &u(0,t)=u(1,t)=0,~~0\leq t\leq 1,\\
 &u(x,0)=0,~~0\leq x\leq 1,
 \end{array}$$
 with  the exact solution given by $u(x,t)=t^2\sin(2\pi x)$.
 \end{example}

 \begin{table}[h]
 \begin{center}
 \caption{Numerical convergence orders in temporal direction with $h=\frac1{100}$ for Example \ref{ex1}.}

 \bigskip
 \label{table1}
 \noindent
 \renewcommand{\arraystretch}{0.9}
 {\small
 \begin{tabular}{lllll}
 \hline
 $\alpha$     & $\tau$   & $E_\infty(h,\tau)$ & $E_2(h,\tau)$  & $Rate1$\\\hline
 $\alpha=0.3$ & 1/5      & 5.9575e-3          & 4.2126e-3      & $\ast$\\
              & 1/10     & 1.4796e-3          & 1.0463e-3      & 2.0095\\
              & 1/20     & 3.6392e-4          & 2.5733e-4      & 2.0235\\
              & 1/40     & 9.1188e-5          & 6.4480e-5      & 1.9967\\
              & 1/80     & 2.2869e-5          & 1.6171e-5      & 1.9954\\
              & 1/160    & 5.7679e-6          & 4.0785e-6      & 1.9873\\[4pt]
 $\alpha=0.5$ & 1/5      & 1.0436e-2          & 7.3790e-3      & $\ast$\\
              & 1/10     & 2.6295e-3          & 1.8594e-3      & 1.9886\\
              & 1/20     & 6.5613e-4          & 4.6395e-4      & 2.0028\\
              & 1/40     & 1.6477e-4          & 1.1651e-4      & 1.9935\\
              & 1/80     & 4.1396e-5          & 2.9272e-5      & 1.9929\\
              & 1/160    & 1.0424e-5          & 7.3708e-6      & 1.9896\\[4pt]
 $\alpha=0.7$ & 1/5      & 1.4850e-2          & 1.0500e-2      & $\ast$\\
              & 1/10     & 3.8141e-3          & 2.6970e-3      & 1.9610\\
              & 1/20     & 9.6191e-4          & 6.8018e-4      & 1.9874\\
              & 1/40     & 2.4234e-4          & 1.7136e-4      & 1.9889\\
              & 1/80     & 6.1157e-5          & 4.3245e-5      & 1.9865\\
              & 1/160    & 1.5448e-5          & 1.0923e-5      & 1.9851\\
 \hline
 \end{tabular}
 }
 \end{center}
 \end{table}

 \begin{table}[h]
 \begin{center}
 \caption{Numerical convergence orders in spatial direction with $\tau=\frac1{4000}$ when $\alpha=0.5$ for Example \ref{ex1}.}

 \bigskip
 \label{table2}
 \noindent
 \renewcommand{\arraystretch}{0.9}
 {\small
 \begin{tabular}{llll}
 \hline
 $h$   & $E_\infty(h,\tau)$ & $E_2(h,\tau)$ & $Rate2$\\\hline
 1/4   & 2.7020e-2          & 1.9106e-2          & $\ast$\\
 1/8   & 1.5651e-3          & 1.1067e-3          & 4.1097\\
 1/16  & 9.6041e-5          & 6.7911e-5          & 4.0265\\
 1/32  & 5.9906e-6          & 4.2360e-6          & 4.0029\\
 1/64  & 3.8961e-7          & 2.7549e-7          & 3.9426\\
 \hline
 \end{tabular}
 }
 \end{center}
 \end{table}

 Convergence order of the proposed scheme in temporal direction with $h=\frac1{100}$ is reported in  Table \ref{table1},
 while in Table \ref{table2}, the convergence order in spatial direction with $\tau=\frac1{4000},~\alpha=0.5$ is listed.
 In this particular example, since the convergent rates of the maximum norm and the 2-norm coincide for the digits listed, they are reported as $Rate1$ and $Rate2$ in the tables. Once again, both tables confirm the theoretical result.

 \bigskip
 Our next two examples are for fractional diffusion-wave equations:
 \begin{example}\label{ex3}
 Consider the problem:
 $$\begin{array}{rl}
 &_{0}^CD_t^{\gamma}u
 =\dfrac{\partial^2u}{\partial x^2}+e^x[\Gamma(\gamma+2)t-t^{\gamma+1}],
  \quad 0\leq x\leq 1, \quad 0<t\leq 1, \quad 1<\gamma<2,\\[5pt]
 &u(x,0)=\dfrac{\partial u(x,0)}{\partial t}=0, \quad 0\leq x\leq 1,\\
 &u(0,t)=t^{1+\gamma}, \quad u(1,t)=et^{1+\gamma}, \quad 0<t\leq 1,
 \end{array}$$
 the problem can be equivalently changed to
 $$\begin{array}{rl}
 &\dfrac{\partial u(x,t)}{\partial t}={_{0}I^\alpha_tu_{xx}(x,t)}
 +e^x\bigg[(\alpha+2)t^{\alpha+1}
 -\frac{\Gamma(\alpha+3)}{\Gamma(2\alpha+3)}t^{2\alpha+2}\bigg],
 \quad 0\leq x\leq 1, \quad 0<t\leq 1,\\
 &u(x,0)=\dfrac{\partial u(x,0)}{\partial t}=0, \quad 0\leq x\leq 1,\\
 &u(0,t)=t^{2+\alpha}, \quad u(1,t)=et^{2+\alpha}, \quad 0<t\leq 1,
 \end{array}$$
 where $\alpha=\gamma-1$. The exact solution for this problem is $u(x,t)=e^xt^{\alpha+2}$.
 \end{example}

 \begin{table}[h]
 \begin{center}
 \caption{Numerical convergence orders in temporal direction with $h=\frac1{30}$ for Example \ref{ex3}.}

 \bigskip
 \label{table7}
 \noindent
 \renewcommand{\arraystretch}{0.9}
 {\small
 \begin{tabular}{llllll}
 \hline
 $\alpha$     & $\tau$   & $E_\infty(h,\tau)$ & $Rate1_\infty$   & $E_2(h,\tau)$ & $Rate1_2$\\\hline
 $\alpha=0.3$ & 1/5      & 5.2645e-3          & $\ast$  & 3.9527e-3     & $\ast$\\
              & 1/10     & 1.3874e-3          & 1.9239  & 9.8966e-4     & 1.9978\\
              & 1/20     & 3.7394e-4          & 1.8916  & 2.6259e-4     & 1.9141\\
              & 1/40     & 9.8549e-5          & 1.9239  & 6.9307e-5     & 1.9218\\
              & 1/80     & 2.5742e-5          & 1.9367  & 1.8173e-5     & 1.9312\\
              & 1/160    & 6.6877e-6          & 1.9445  & 4.7249e-6     & 1.9435\\[4pt]
 $\alpha=0.6$ & 1/5      & 8.6638e-3          & $\ast$  & 5.9115e-3     & $\ast$\\
              & 1/10     & 2.2702e-3          & 1.9322  & 1.5390e-3     & 1.9415\\
              & 1/20     & 5.7192e-4          & 1.9890  & 3.8707e-4     & 1.9914\\
              & 1/40     & 1.4463e-4          & 1.9835  & 9.7399e-5     & 1.9906\\
              & 1/80     & 3.6248e-5          & 1.9963  & 2.4487e-5     & 1.9919\\
              & 1/160    & 9.0979e-6          & 1.9943  & 6.1472e-6     & 1.9940\\[4pt]
 $\alpha=0.9$ & 1/5      & 1.3142e-2          & $\ast$  & 8.7666e-3     & $\ast$\\
              & 1/10     & 3.4817e-3          & 1.9163  & 2.2332e-3     & 1.9729\\
              & 1/20     & 8.8458e-4          & 1.9767  & 5.5942e-4     & 1.9971\\
              & 1/40     & 2.2304e-4          & 1.9877  & 1.3978e-4     & 2.0008\\
              & 1/80     & 5.5750e-5          & 2.0003  & 3.4943e-5     & 2.0001\\
              & 1/160    & 1.3919e-5          & 2.0019  & 8.7372e-5     & 1.9998\\
 \hline
 \end{tabular}
 }
 \end{center}
 \end{table}

 \begin{table}[ht]
 \begin{center}
 \caption{Numerical convergence orders in spatial direction with $\tau=\frac1{5000}$ when $\alpha=0.5$ for Example \ref{ex3}.}

 \bigskip
 \label{table8}
 \noindent
 \renewcommand{\arraystretch}{0.9}
 {\small
 \begin{tabular}{lllll}
 \hline
 $h$   & $E_\infty(h,\tau)$ & $Rate2_\infty$& $E_2(h,\tau)$ & $Rate2_2$\\\hline
 1/2   & 3.6421e-5          & $\ast$        & 2.5753e-5     & $\ast$\\
 1/4   & 2.2865e-6          & 3.9936        & 1.7004e-6     & 3.9208\\
 1/8   & 1.4229e-7          & 4.0062        & 1.0616e-7     & 4.0016\\
 1/16  & 8.5477e-9          & 4.0572        & 5.9751e-9     & 4.1511\\
 \hline
 \end{tabular}
 }
 \end{center}
 \end{table}

 Table \ref{table7} and Table \ref{table8} justify the accuracy of the scheme proposed in Section 4.

 \bigskip
 In the last example,
 we consider a problem where the exact solution cannot be found readily. We note that, in this example, $u_t$ is not identically equal to zero initially.
 \begin{example}\label{ex4}
 Consider the problem:
 $$\begin{array}{rl}
 &_{0}^CD_t^{\gamma}u
 =\dfrac{\partial^2u}{\partial x^2}+\sin(2\pi x)
 \bigg[\frac{\Gamma(\gamma+3)}2t^2+4\pi^2t^{2+\gamma}\bigg],
 \quad 0\leq x\leq 1, \quad 0<t\leq 1, \quad 1<\gamma<2,\\[5pt]
 &u(x,0)=0,\quad \dfrac{\partial u(x,0)}{\partial t}=0.1\sin(2\pi x), \quad 0\leq x\leq 1,\\
 &u(0,t)=u(1,t)=0, \quad 0<t\leq 1.
 \end{array}$$
 By noting $\alpha=\gamma-1$, the problem can be equivalently changed to
 $$\begin{array}{l}
 \dfrac{\partial u(x,t)}{\partial t}={_{0}I^\alpha_tu_{xx}(x,t)}
 +\sin(2\pi x)\bigg[(\alpha+3)t^{\alpha+2}
 +4\pi^2\frac{\Gamma(\alpha+4)}{\Gamma(2\alpha+4)}t^{2\alpha+3}\bigg],
 \quad 0\leq x\leq 1, \quad 0<t\leq 1,\\
 u(x,0)=0,\quad \dfrac{\partial u(x,0)}{\partial t}=0.1\sin(2\pi x), \quad 0\leq x\leq 1,\\
 u(0,t)=u(1,t)=0, \quad 0<t\leq 1.
 \end{array}$$
 \end{example}

 We take the numerical solution with $M=400$, $N=4000$ as the `true' solution when computing the errors.
 From the last two tables,
 we can see that the scheme still works properly in this situation.

 \begin{table}[h]
 \begin{center}
 \caption{Numerical convergence orders in temporal direction with $h=\frac1{50}$ for Example \ref{ex4}.}

 \bigskip
 \label{table-ini}
 \noindent
 \renewcommand{\arraystretch}{0.9}
 {\small
 \begin{tabular}{llllll}
 \hline
 $\alpha$  & $\tau$  & $E_\infty(h,\tau)$ & $E_2(h,\tau)$ & $Rate1$\\\hline
 $\alpha=0.3$ & 1/10   & 2.9448e-3   & 2.0864e-3  & $\ast$\\
              & 1/20   & 8.1395e-4   & 5.7669e-4  & 1.8551\\
              & 1/40   & 2.1283e-4   & 1.5079e-4  & 1.9353\\
              & 1/80   & 5.3775e-5   & 3.8100e-5  & 1.9847\\
              & 1/160  & 1.2898e-5   & 9.1384e-6  & 2.0598\\[4pt]
 $\alpha=0.5$ & 1/10   & 4.6281e-3   & 3.2791e-3  & $\ast$\\
              & 1/20   & 1.3211e-3   & 9.3602e-4  & 1.8087\\
              & 1/40   & 3.4933e-4   & 2.4750e-4  & 1.9191\\
              & 1/80   & 8.9065e-5   & 6.3103e-5  & 1.9717\\
              & 1/160  & 2.1877e-5   & 1.5500e-5  & 2.0254\\[4pt]
 $\alpha=0.7$ & 1/10   & 6.4028e-3   & 4.5364e-3  & $\ast$\\
              & 1/20   & 1.6342e-3   & 1.1579e-3  & 1.9701\\
              & 1/40   & 4.2105e-4   & 2.9832e-4  & 1.9565\\
              & 1/80   & 1.0677e-4   & 7.5650e-5  & 1.9794\\
              & 1/160  & 2.6331e-5   & 1.8655e-5  & 2.0197\\
 \hline
 \end{tabular}
 }
 \end{center}
 \end{table}

 \begin{table}[ht]
 \begin{center}
 \caption{Numerical convergence orders in spatial direction with $\tau=\frac1{2000}$ when $\alpha=0.5$ for Example \ref{ex4}.}

 \bigskip
 \label{table-ini2}
 \noindent
 \renewcommand{\arraystretch}{0.9}
 {\small
 \begin{tabular}{lllll}
 \hline
 $h$   & $E_\infty(h,\tau)$ & $Rate2_\infty$& $E_2(h,\tau)$ & $Rate2_2$\\\hline
 1/5    & 9.0588e-3          & $\ast$        & 6.7352e-3     & $\ast$\\
 1/10   & 5.4015e-4          & 4.0679        & 4.0160e-4     & 4.0679\\
 1/20   & 3.4978e-5          & 3.9488        & 2.4733e-5     & 4.0212\\
 1/40   & 2.0755e-6          & 4.0749        & 1.4676e-6     & 4.0749\\
 \hline
 \end{tabular}
 }
 \end{center}
 \end{table}


\begin{thebibliography}{99}

 \bibitem{Podlubny} I. Podlubny, Fractional Differential Equations, Academic Press, New York, 1999.

 \bibitem{Kilbas} A. Kilbas, H. Srivastava, J. Trujillo, Theory and Applications of Fractional Differential Equations, Elsevier Science and Technology, 2006.

 \bibitem{Klages} R. Klages, G. Radons, I. Sokolov, Anomalous Transport: Foundations and Applications, Wiley-VCH, Weinheim, 2008.

 \bibitem{Nigmatullin1} R. Nigmatullin, To the theoretical explanation of the universal response, Phys. Status Solidi (B): Basic Res. 123 (2) (1984) 739--745.

 \bibitem{Nigmatullin2} R. Nigmatullin, Realization of the generalized transfer equation in a medium with fractal geometry, Phys. Status Solidi (B): Basic Res. 133 (1) (1986) 425--430.

 \bibitem{ChenC} C. Chen, F. Liu, I. Turner, V. Anh, A Fourier method for the fractional diffusion equation describing sub-diffusion, J. Comput. Phys. 227 (2007) 886--897.

 \bibitem{ZhaoJ} J. Zhao, T. Zhang, R. Corless, Convergence of the compact finite difference method for second-order elliptic equations, Appl. Math. Comput. 182 (2006) 1454--1469.

 \bibitem{SunWu} Z. Sun, X. Wu, A fully discrete difference scheme for a diffusion-wave system, Appl. Numer. Math. 56 (2006) 193--209.

\bibitem{Du} R. Du, W. Cao, Z. Sun, A compact difference scheme for the fractional diffusion-wave equation, Appl. Math. Model. 34 (2010) 2998--3007.

 \bibitem{LiaoW} W. Liao, J. Zhu, A.Q.M. Khaliq, An efficient high-order algorithm for solving systems of reaction-diffusion equations, Numer. Meth. Partial Diff. Eq. 18 (2002) 340--354.

 \bibitem{ZhuangP} P. Zhuang, F. Liu, V. Anh, I. Turner, New solution and analytical techniques of the implicit numerical method for the anomalous subdiffusion equation, SIAM J. Numer. Anal. 46 (2) (2008) 1079--1095.

 \bibitem{r1} Z. Sun, Compact difference schemes for heat equation with Neumann boundary conditions, Numer. Meth. Part. Differ. Equ. 25 (2009) 1320--1341.

 \bibitem{r2} C. Chen, F. Liu, V. Anh, I. Turner, Numerical schemes with high spatial accuracy for a variable-order anomalous subdiffusion equation, SIAM J. Sci. Comput. 32 (2010) 1740--1760.

 \bibitem{r3} C. Chen, F. Liu, I. Turner, V. Anh, Numerical methods with fourth-order spatial accuracy for variable-order nonlinear Stokes' first problem for a heated generalized second grade fluid, Comput. Math. Appl. 62 (2011) 971--986.

 \bibitem{r4} C. Chen, F. Liu, V. Anh, I. Turner, Numerical methods for solving a two-dimensional variable-order anomalous subdiffusion equation, Math. Comput. 82 (2012) 345--366.

 \bibitem{r5} M. Cui, Compact alternating direction implicit method for two-dimensional time fractional diffusion equation, J. Comput. Phys. 231 (2012) 2621--2633.

 \bibitem{r6} J. Ren, Z. Sun, Numerical algorithm with high spatial accuracy for the fractional diffusion-wave equation with Neumann boundary conditions, J. Sci. Comput. 56 (2013) 381--408.

 \bibitem{Langlands} T. Langlands, B. Henry, The accuracy and stability of an implicit solution method for the fractional diffusion equation, J. Comput. Phys. 205 (2005) 719--736.

 \bibitem{Sousaa} E. Sousa, C. Li, A weighted finite difference method for the fractional diffusion equation based on the Riemann-Liouville derivative, (2011), arXiv:1109.2345 [math.NA].

 \bibitem{Meerschaert} M. Meerschaert, C. Tadjeran, Finite difference approximations for fractional advection-dispersion flow equations, J. Comput. Appl. Math. 172 (2004) 65--77.

 \bibitem{Yuste1} S. Yuste, L. Acedo, An explicit finite difference method and a new von Neumann-type stability analysis for fractional diffusion equations, SIAM J. Numer. Anal. 42 (5) (2005) 1862--1874.

 \bibitem{Yuste2}S. Yuste, Weighted average finite difference methods for fractional diffusion equations, J. Comput. Phys. 216 (2006) 264--274.

  \bibitem{Tian} W. Tian, H. Zhou, W. Deng,  A class of second order difference approximations for solving space fractional diffusion equations, Math. Comput. arXiv:1201.5949 [math.NA].

 \bibitem{Zhou} H. Zhou, W. Tian, W. Deng, Quasi-compact finite difference schemes for space fractional diffusion equations, J. Sci. Comput. 56 (2013) 45--66.

 \bibitem{F-Liu} F. Liu, C. Yang, K. Burrage, Numerical method and analytical technique of the modified anomalous subdiffusion equation with a nonlinear source term, J. Comput. Appl. Math. 231 (2009) 160--176.

  \bibitem{Q-Liu} Q. Liu, F. Liu, I. Turner, V. Anh, Finite element approximation for a modified anomalous subdiffusion equation, Appl. Math. Modell., 35 (2011) 4103--4116.

 \bibitem{Mohebbi} A. Mohebbi, M. Abbaszade, M. Dehghan, A high-order and unconditionally stable scheme for the modified anomalous fractional sub-diffusion equation with a nonlinear source term, J. Comput. Phys. 240 (2013) 36--48.

 \bibitem{Lubich} C. Lubich, Discretized fractional calculus, SIAM J. Math. Anal. 17 (1986) 704--719.

 \bibitem{Ervin} V. Ervin, J. Roop, Variational formulation for the stationary fractional advection dispersion equation, Numer. Meth. Part. Differ. Equ. 22 (2006) 558--576.

  \bibitem{Gao} G. Gao, Z. Sun, A compact finite difference scheme for the fractional sub-diffusion equations, J. Comput. Phys. 230 (2011) 586--595.

 \bibitem{Chan2} R. Chan, X. Jin, An Introduction to Iterative Toeplitz Solvers, SIAM (2007).

 \bibitem{Quarteroni} A. Quarteroni, A. Valli, Numerical approximation of partial differential equations. Springer, Berlin, 1997.

 \bibitem{Huang} J. Huang, Y. Tang, L. V\'{a}zquez, J. Yang, Two finite difference schemes for time fractional diffusion-wave equation, Numer. Algor. 64 (2013) 707--720.

 \bibitem{Li} C. Li, A. Chen, J. Ye, Numerical approaches to fractional calculus and fractional ordinary differential equation, J. Comput. Phys. 230 (2011) 3352--3368.

 \bibitem{Cui} M. Cui, Compact finite difference method for the fractional diffusion equation, J. Comput. Phys. 228 (2009) 7792--7804.

 \bibitem{Lin} Y. Lin, C. Xu, Finite difference/spectral approximations for the time-fractional diffusion equation, J. Comput. Phys. 225 (2007) 1533--1552.
\end{thebibliography}
\end{document}